# The distribution of natural numbers divisible by 2, 3, 5, 11, 13 and 17 on the Square Root Spiral

-  Harry K. Hahn  -

Ettlingen / Germany

28. June  2007

**Abstract :**

The natural numbers divisible by the **Prime Factors 2, 3, 5, 11, 13** and **17** lie on defined spiral-graphs, which run through the Square Root Spiral.  And a mathematical analysis shows that these spiral graphs are defined by specific quadratic polynomials.

Basically all natural number which are divisible by the same prime factor lie on such spiral-graphs. And these spiral-graphs can be assigned to a certain number of "spiral graph systems"  which have a defined spatial orientation to each other.

This document represents a supplementation to my detailed introduction study to the Square Root Spiral, and it contains the missing diagrams and analyses for the distribution of the natural numbers divisible by 2, 3, 5, 11, 13 and 17 on the Square Root Spiral.   My introduction study can be found in the ArXiv–archive :

Title of the study :    " **The ordered distribution of the natural numbers on the Square Root Spiral** "

( This study also contains diagram & analysis for numbers divisible by 7  → not included in this study here ! )

**Introduction :**

The Square Root Spiral ( or "Spiral of Theodorus" or "Einstein-Spiral" ) is a very interesting geometrical structure, in which the square roots of all natural numbers have a defined (spatial ) position to each other.

The Square Root Spiral develops from a right angled base triangle with the two legs ( cathets ) having the length 1, and with the long side ( hypotenuse ) having a length which is equal to the square root of 2.
The Square Root Spiral is formed by further adding right angled triangles on this base triangle.
In this process the longer legs of the next triangles always attach to the hypotenuses of the previous triangles. And the longer leg of the next triangle always has the same length as the hypotenuse of the previous triangle, and the shorter leg always has the length 1.

In this way a spiral structure is developing in which the spiral is created by the shorter legs of the triangles, which have the constant length of 1, and where the lengths of the radial rays ( or spokes ) coming from the centre of this spiral are the square roots of the natural numbers ( sqrt 2 , sqrt 3, sqrt 4, sqrt 5 …. ).

The most striking property of the Square Root Spiral is surely the fact, that the distance between two successive winds of the Square Root Spiral quickly strives for the well known geometrical constant  $\pi$   !!

But there are many more interesting interdependencies between the natural numbers, which can be discovered in this amazing geometrical structure.  For example all numbers divisible by the same prime factor always lie on defined spiral-graphs, which run in a clear order through the Square Root Spiral. And the square numbers  4, 9, 16, 25, ...  form a highly three-symmetrical sytem of three spiral-graphs, which divide the Square Root Spiral into three equal areas.

It is the same with Prime Numbers.   Prime Numbers also clearly accumulate on defined spiral-graphs.

A mathematical analysis shows, that all these spiral graphs represent special quadratic polynomials.

Even the famous Fibonacci Numbers play a part in the structure of the Square Root Spiral.  Fibonacci - Numbers divide the Square Root Spiral into areas and angles with constant proportions, which are linked to the "golden mean" ( or golden section ).
The Square Root Spiral ( or Einstein Spiral ) shows the interdependencies between the natural numbers in a visual way.  Therefore it can be considered to be a kind of visual representation of number theory !

Through pure graphical analysis of this amazing structure, the higher logic of the ( spatial ) distribution of the natural numbers ( and special sub-groups like square numbers or prime numbers ) comes to light and is very easy to grasp, because it is visible !

Please have a read through my more detailed introduction about the distribution of the natural numbers on the Square Root Spiral.  This introduction gives a good insight into the interdependencies between the natural numbers,  which exist in this amazing geometrical structure.

→  Titel :     "The ordered distribution of the natural numbers on the Square Root Spiral"

## The distribution of numbers divisible by 2 on the square root spiral :

If I am talking about the arrangement of numbers divisible by **2** on the square root spiral, I am actually referring to the square roots of these numbers, which appear as radial rays on the square root spiral. The natural numbers divisible by **2** themselves, can be seen as imaginary square areas staying vertical on these radial rays.

( → explanation see in my detailed introduction study :   " The ordered distribution of the natural numbers on the square root spiral" )

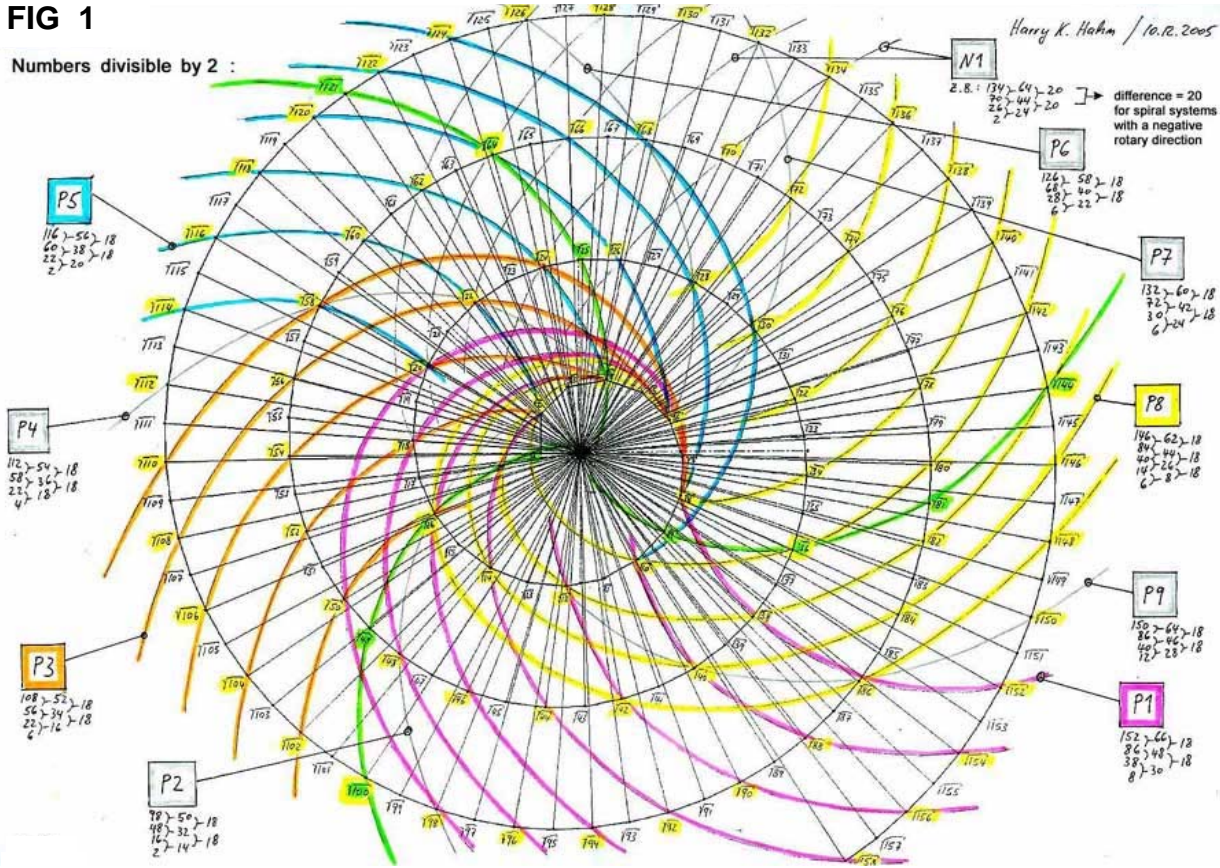

**FIG 1**

From the image **FIG 1** it is evident that the numbers divisible by **2** ( marked in yellow ) lie on defined spiral graphs which have their starting point in or near the centre of the square root spiral.

These spiral graphs have either a positive or a negative rotation direction.

A spiral graph which has a clockwise rotation direction shall be called negative (N) and a spiral graph which has an counterclockwise rotation direction shall be called positive (P).

The green spiral-graphs show the highly three-symmetrical arrangement of the three spiral-graphs which contain the square numbers 4, 9, 16, 25, 36, 49, …. ( drawn for reference only ! )

**Here now a listing of the most important properties of the found spiral graphs :**

- **The numbers divisible by 2** lie on **10** spiral-graph-systems with a **negative** rotation direction and on **9** spiral-graph-systems with a **positive** rotation direction.
  From the **10** spiral-graph-systems with a negative rotation direction only the system **N1** is drawn in light grey color. And from the **10** spiral graph systems with a positive rotation direction only the systems **P1, P3, P5** and **P8** are drawn in more detail ( drawn in the colors pink, orange, blue and yellow ).  From the other spiral graph systems with a positive rotation direction **P2, P4, P6, P7** and **P9** only one spiral arm for each system is drawn in favour of clearness. These spiral arms are drawn in light grey and they lie between the other four spiral graph systems which are drawn in color.

- The **10** negative spiral graph systems **N1 – N10** are arranged in an angle of around 36° to each other, whereas the **9** positive spiral graph systems **P1 – P9** are arranged in an angle of around 360°/9 = 40° to each other ( in reference to the centre of the square root spiral )

    From the negative Spiral Graph Systems the systems **N1** and **N6, N2** and **N7, N3** and **N8, N4** and **N9** as well as **N5** and **N10,** form pairs of spiral graph systems which lie approximately point-symmetrical to each other ( in reference to the centre of the square root spiral ).
    However the positive spiral graph systems **P1 – P9** don't form such point-symmetrical pairs.

- Calculating the differences of the successive numbers lying on one of the drawn spiral arms, and then further calculating the differences of these differences, results in the constant value of **18** for the spiral graph systems with a positive rotation direction ( **P1 – P9** ), and it results in the constant value of **20** for the spiral graph systems with a negative rotation direction( **N1 – N10** ).

    ( The result of this repeated difference calculation shall be called **2. Differential** )

    → see difference calculation for 10 exemplary spiral arms ( one spiral arm per system )
    in **FIG 1** beside the names of the spiral graph systems P1 to P9 and N1.
    ( → systems  N2 to N10  not shown ! )

    These 10 exemplary spiral arms are defined by the following quadratic polynomials :

    **Quadratic Polynomials of exemplary spiral arms :**

    | | | |
    |---|---|---|
    | $9x^2 + 21x + 8$ | → | belongs to **P1** - system |
    | $9x^2 + 5x + 2$ | → | belongs to **P2** - system |
    | $9x^2 + 7x + 6$ | → | belongs to **P3** - system |
    | $9x^2 + 9x + 4$ | → | belongs to **P4** - system |
    | $9x^2 + 11x + 2$ | → | belongs to **P5** - system |
    | $9x^2 + 13x + 6$ | → | belongs to **P6** - system |
    | $9x^2 + 15x + 6$ | → | belongs to **P7** - system |
    | $9x^2 + 17x + 14$ | → | belongs to **P8** - system |
    | $9x^2 + 19x + 12$ | → | belongs to **P9** - system |
    | $10x^2 + 14x + 2$ | → | belongs to **N1** - system |

**Note :**  The positive (P) and the negative (N) spiral graph systems of the natural numbers, divisible by 2 , have a  **2. Differential** which is not identical !

**The distribution of numbers divisible by 3 on the square root spiral :**

## FIG 2

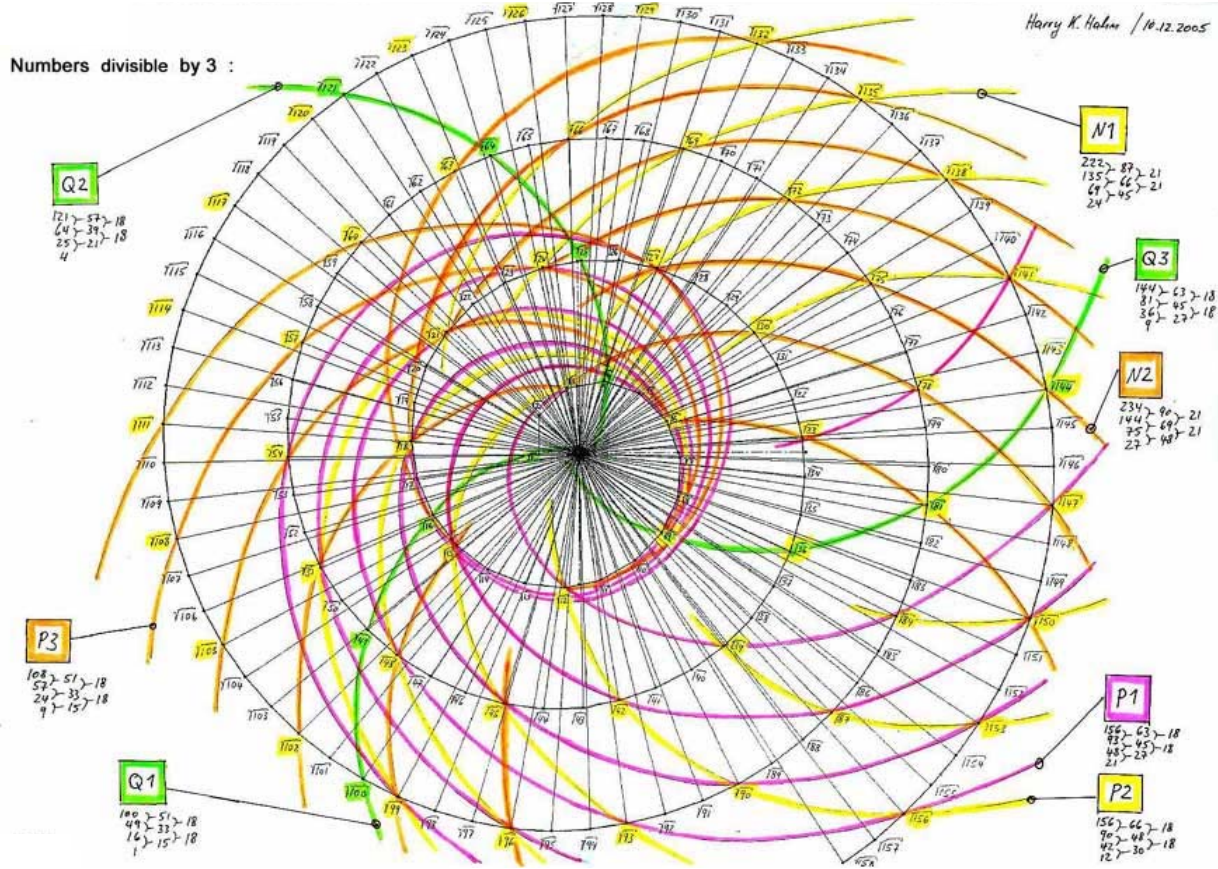

The numbers divisible by **3** ( marked in yellow ) lie on defined spiral graphs which have their starting point in or near the centre of the square root spiral.

These spiral graphs have either a positive or a negative rotation direction. A spiral graph which has a clockwise rotation direction shall be called negative (N) and a spiral graph which has an anticlockwise rotation direction shall be called positive (P).

The green spiral-graphs show the highly three-symmetrical arrangement of the three spiral-graphs which contain the square numbers 4, 9, 16, 25, 36, 49, …. ( drawn for reference only ! )

**Here now a listing of the most important properties of the found spiral graphs :**

- **The numbers divisible by 3** lie on **7** spiral-graph-systems with a **negative** rotation direction ( only two of these spiral graph systems are drawn in orange and yellow ) and on **6** spiral graph systems with a **positive** rotation direction ( **Note :** only three of these spiral graph systems are drawn in pink, yellow and orange ).
  The **7** negative spiral graph systems are named **N1 – N7** and the **6** positive systems are named **P1 – P6** ( only the two negative spiral-systems **N1** and **N2** and the three positive spiral systems **P1** to **P3** are drawn ! ).

- The **6** positive spiral graph systems **P1 – P6** are arranged in an angle of around 60° to each other, whereas the **7** negative spiral graph systems **N1 – N7** are arranged in an angle of around 360°/7 = 51.43° to each other ( in reference to the centre of the Square Root Spiral ).
  From the positive spiral graph systems, the systems **P1** and **P4**, **P2** and **P5** as well as **P3** and **P6** form pairs of spiral graph systems which lie approximately point-symmetrical to each other ( in reference to the centre of the square root spiral ).
  However the negative spiral graph systems **N1 – N7** don't form such point-symmetrical pairs.

- Calculating the differences of the successive numbers lying on one of the drawn spiral arms, and then further calculating the differences of these differences, results in the constant value of **18** for the spiral graph systems with a positive rotation direction ( **P1 – P6** ), and it results in the constant value of **21** for the spiral graph systems with a negative rotation direction( **N1 – N7** ).

   ( The result of this repeated difference calculation shall be called **2. Differential** )

   → **see difference calculation for 5 exemplary spiral arms** ( one spiral arm per system )
      **in FIG 2** beside the names of the spiral graph systems N1, N2 and P1, P2, P3 .
      ( → systems N3 to N7 and P4 to P6 not shown ! )

   These 5 exemplary spiral arms are defined by the following quadratic polynomials :

   **Quadratic Polynomials of exemplary spiral arms :**

   | | | |
   |---|---|---|
   | $10.5x^2 + 34.5x + 24$ | → | belongs to **N1** - system |
   | $10.5x^2 + 37.5x + 27$ | → | belongs to **N2** - system |
   | $9x^2 + 12$ | → | belongs to **P1** - system |
   | $9x^2 + 21x + 12$ | → | belongs to **P2** - system |
   | $9x^2 + 6x + 9$ | → | belongs to **P3** - system |

**Note :**  The positive (P) and the negative (N) spiral graph systems of the natural numbers, divisible by 3 , have a **2. Differential** which is not identical !

**The distribution of numbers divisible by 5 on the square root spiral :**

**FIG 3**

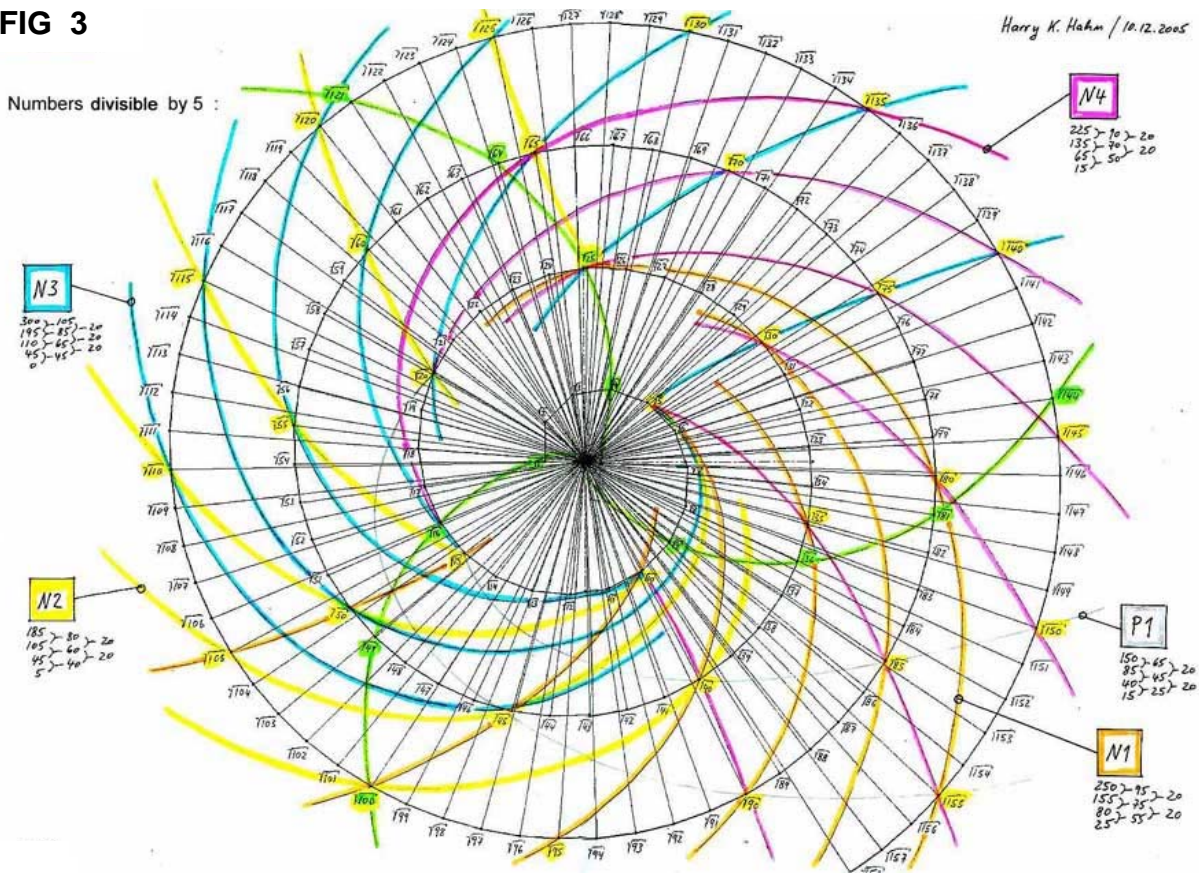

The numbers divisible by **5** ( marked in yellow ) lie on defined spiral graphs which have their starting point in or near the centre of the square root spiral.

The green spiral-graphs show the highly three-symmetrical arrangement of the three spiral-graphs which contain the square numbers 4, 9, 16, 25, 36, 49, …. ( drawn for reference only ! )

**Here now a listing of the most important properties of the found spiral graphs :**

- **The numbers divisible by 5** lie on **4** spiral graph systems with a **negative** rotation direction ( drawn in orange, yellow, blue and pink ) and on **4** spiral graph systems with a **positive** rotation direction ( only one system ( P1 ) drawn in light gray lines ! ).
  The **4** negative spiral graph systems are named **N1, N2, N3** and **N4** and the **4** positive systems are named **P1, P2, P3** and **P4** ( only two spiral arms of P1 are drawn ! ).

- The spiral graph systems **N1, N2, N3** and **N4** as well as **P1, P2, P3** and **P4** are arranged in an angle of around 90° to each other. And two at a time of these spiral graph systems are located to each other in a point-symmetrical manner ( in reference to the centre of the square root spiral ).
  For example the two systems **N1** and **N3** , as well as the systems **N2** and **N4** seem to lie point-symmetrical to each other.
  The same applies to the spiral graph systems **P1 / P3** and **P2 / P4** ( only P1 drawn ! )

- Calculating the differences of the successive numbers lying on one of the drawn spiral arms, and then further calculating the differences of these differences, results in the constant value of **20** for the positive as well as the negative rotating spiral-graphs. ( **2. Differential = 20** )

    → **see difference calculation for 5 exemplary spiral arms** ( one spiral arm per system )
      **in FIG 3** beside the names of the spiral graph systems N1, N2, N3, N4 and P1.
      ( → systems P2 to P4 not shown ! )

    These 5 exemplary spiral arms are defined by the following quadratic polynomials :

**Quadratic Polynomials of exemplary spiral arms :**

| | |
|---|---|
| **$10x^2 + 45x + 25$**   → | belongs to **N1** - system |
| **$10x^2 + 30x + 5$**   → | belongs to **N2** - system |
| **$10x^2 + 55x + 45$**   → | belongs to **N3** - system |
| **$10x^2 + 40x + 15$**   → | belongs to **N4** - system |
| **$10x^2 + 15x + 15$**   → | belongs to **P1** - system |

# The distribution of numbers divisible by 11 on the Square Root Spiral

**FIG 4**

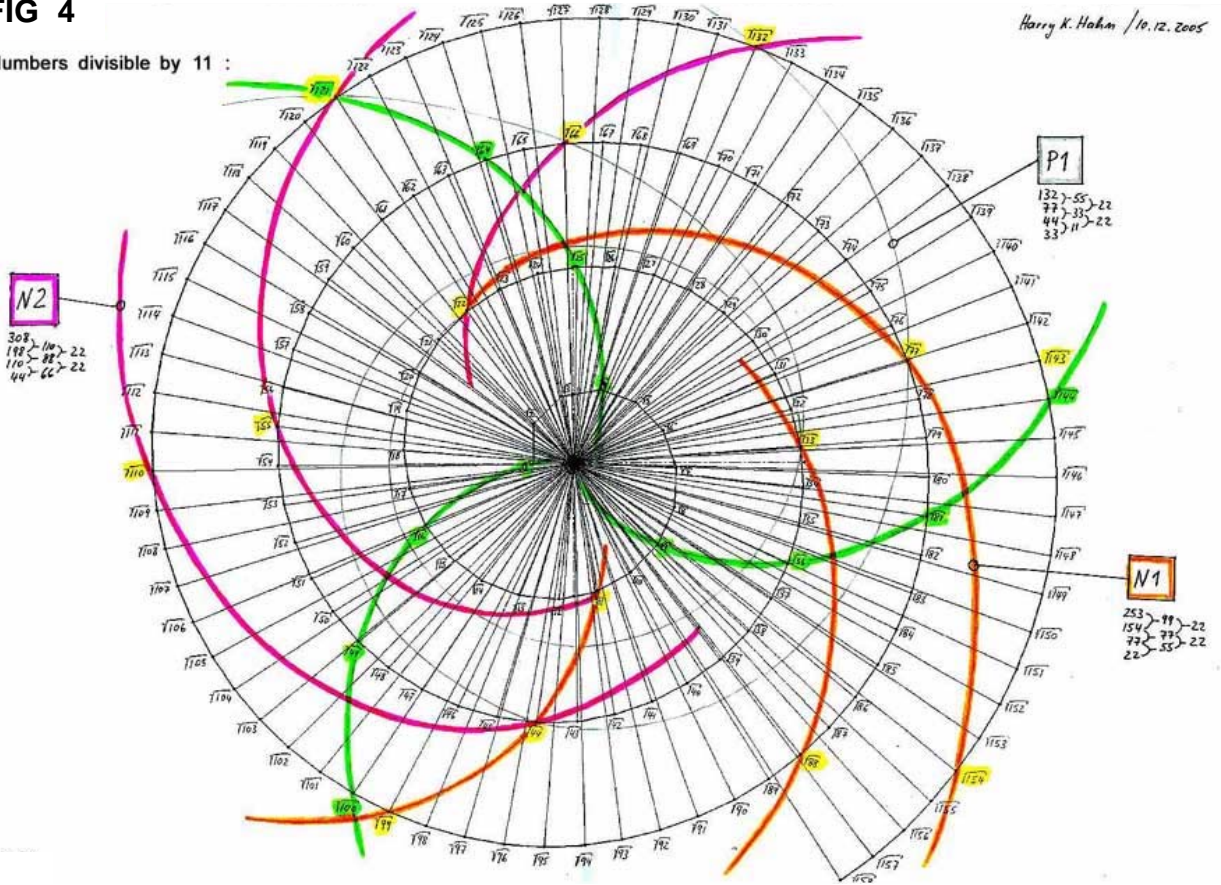

The numbers divisible by **11** ( marked in yellow ) lie on defined spiral graphs which start in or near the centre of the square root spiral. These spiral graphs have either a positive or a negative rotation direction.

The green spiral-graphs show the highly three-symmetrical arrangement of the three spiral-graphs which contain the square numbers 4, 9, 16, 25, 36, 49, …. ( drawn for reference only ! )

**Here now a listing of the most important properties of the found spiral graphs :**

- **The numbers divisible by 11** lie on **2** spiral graph systems with a **negative** rotation direction ( drawn in orange and pink ) and on **2** spiral graph systems with a **positive** rotation direction ( only one system with positive rotation direction drawn in light gray lines ! ). The **2** negative spiral graph systems are named **N1** and **N2** and the **2** positive systems are named **P1** and **P2**. ( → only P1 is drawn ! ).
  
  Note :   Not all spiral arms of the described spiral graph systems are drawn !

- The spiral graph systems **N1** and **N2** as well as **P1** and **P2** lie approximately point-symmetrical to each other ( in reference to the centre of the square root spiral )

- Calculating the differences of the successive numbers lying on one of the drawn spiral arms, and then further calculating the differences of these differences, results in the constant value of **22** for the positive as well as the negative rotating spiral-graphs.
  
  ( The result of this repeated difference calculation shall be called **2. Differential** )
  
  → see difference calculation for 3 exemplary spiral arms ( one spiral arm per system )
     in FIG 2 beside the names of the spiral graph systems N1, N2 and P1.
     ( → P2-system not shown ! )
  
  These 3 exemplary spiral arms are defined by the following quadratic polynomials :

**Quadratic Polynomials of exemplary spiral arms :**

| | | |
|---|---|---|
| $11x^2 + 44x + 22$ | → | belongs to **N1** – system |
| $11x^2 + 33x + 11$ | → | belongs to **N2** – system |
| $11x^2 + 22x + 44$ | → | belongs to **P1** - system |

# The distribution of numbers divisible by 13 on the Square Root Spiral

**FIG 5**

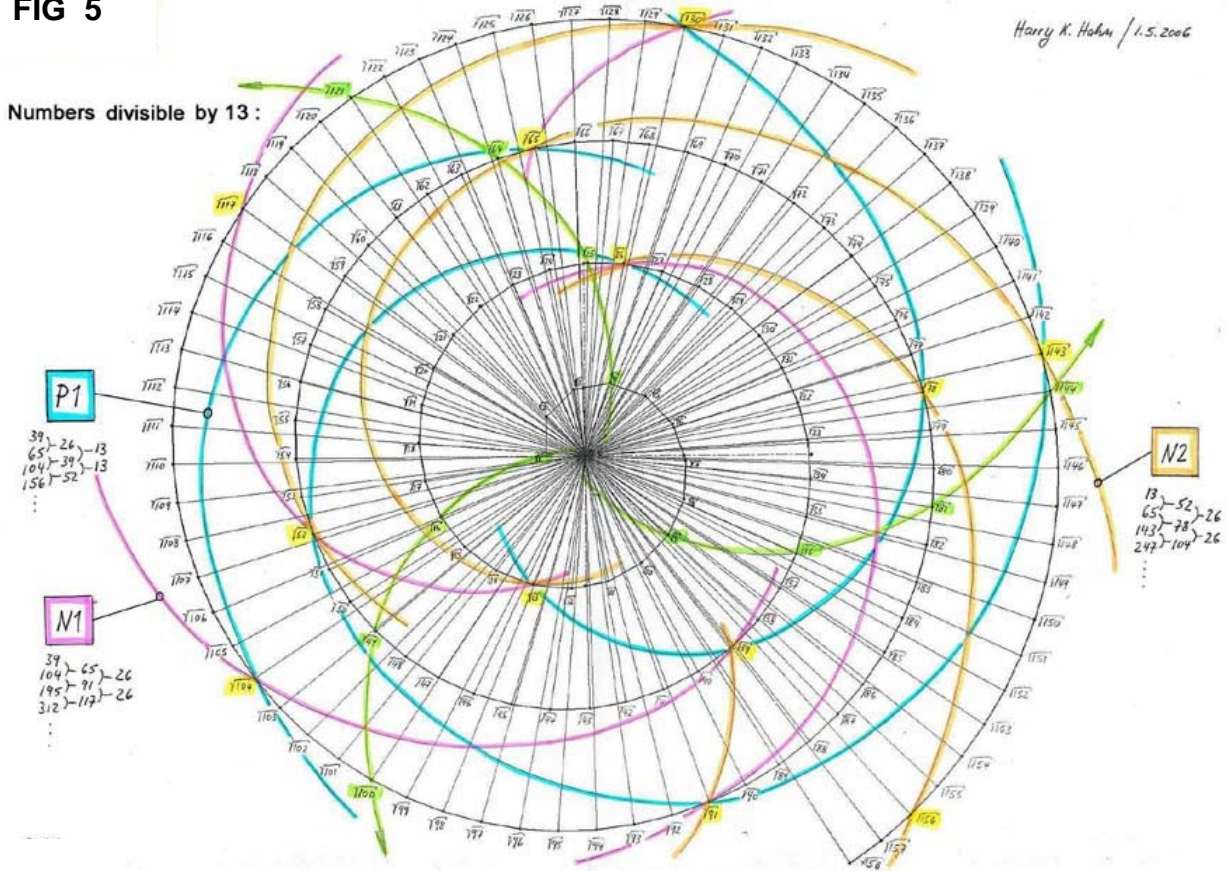

Numbers divisible by 13 :

The numbers divisible by **13** ( marked in yellow ) lie on defined spiral graphs which start in or near the centre of the square root spiral. These spiral-graphs have either a positive or a negative rotation direction.

The green spiral-graphs show the highly three-symmetrical arrangement of the three spiral-graphs which contain the square numbers 4, 9, 16, 25, 36, 49, …. ( drawn for reference only ! )

**Here now a listing of the most important properties of the found spiral graphs :**

- **The numbers divisible by 13** lie on **2** spiral graph systems with a **negative** rotation direction ( drawn in pink and orange ) and on **1** spiral graph system with a **positive** rotation direction ( drawn in blue ). The **2** negative spiral graph systems are named **N1** and **N2** and the positive system is named **P1**.

- The spiral graph systems **N1** and **N2** lie **axis-symmetrical** to each other ( in reference to an imaginary axis running approximately from sqrt 116 to sqrt 152 )

- Calculating the differences of the successive numbers lying on one of the drawn spiral arms, and then further calculating the differences of these differences, results in the constant value of **13** for the spiral graph system with the **positive** rotation direction ( **P1** ), and it results in the constant value of **26** for the two spiral graph systems with a **negative** rotation direction ( **N1** and **N2** ).

  ( The result of this repeated difference calculation shall be called **2. Differential** )

  → see difference calculation for 3 exemplary spiral arms ( one spiral arm per system )
     in FIG 3 beside the names of the spiral graph systems N1, N2 and P1.

  These 3 exemplary spiral arms are defined by the following quadratic polynomials :

  **Quadratic Polynomials of exemplary spiral arms :**

  $13x^2 + 52x + 39$ → belongs to **N1** – system
  $13x^2 + 39x + 13$ → belongs to **N2** – system
  $6.5x^2 + 6.5x + 26$ → belongs to **P1** - system

# The distribution of numbers divisible by 17 on the Square Root Spiral

**FIG 6**

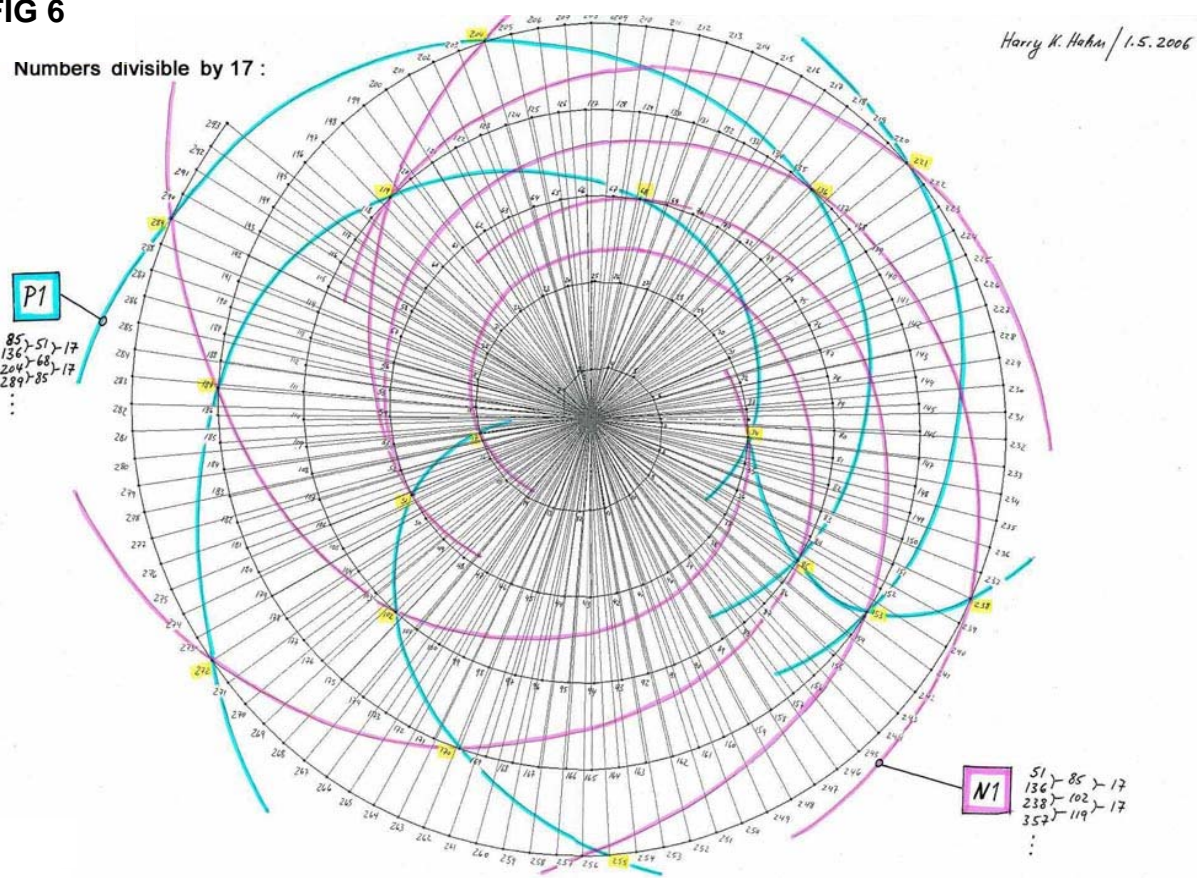

The numbers divisible by **17** ( marked in yellow ) lie on defined spiral graphs which have their starting point in or near the centre of the square root spiral. These spiral graphs have either a positive or a negative rotation direction.

A spiral graph which has a clockwise rotation direction shall be called negative (N) and a spiral graph which has an anticlockwise rotation direction shall be called positive (P).

**Here now a listing of the most important properties of the found spiral graphs :**

- **The numbers divisible by 17** lie on **1** spiral graph system with a **negative** rotation direction ( drawn in pink ) and on **1** spiral graph system with a **positive** rotation direction ( drawn in blue ).

- Calculating the differences of the successive numbers lying on one of the drawn spiral arms, and then further calculating the differences of these differences, results in the constant value of **17** for the spiral graph system with the **positive** rotation direction ( **P1** ), and it results in the same constant value of **17** for the spiral graph system with the **negative** rotation direction ( **N1** ).

  ( → **2. Differential = 17** )

  → see difference calculation for 2 exemplary spiral arms ( one spiral arm per system )
  in **FIG 4** beside the names of the spiral graph systems N1 and P1.

  These 2 exemplary spiral arms are defined by the following quadratic polynomials :

  **Quadratic Polynomials of exemplary spiral arms :**

  $8.5x^2 + 76.5x + 51$         →         belongs to **N1** - system

  $8.5x^2 + 8.5x + 34$          →         belongs to **P1** - system

Please also have a read through my more detailed introduction study about the distribution of the natural numbers on the Square Root Spiral. This introduction gives an overall view of the "Spiral-Graphs" belonging to the numbers divisible by the prime factors 2, 3, 5, 7, 11, 13, and 17, as well as an mathematical explanation of the underlying rules which determine the existence of these spiral-graphs.

→   Titel :        **"The ordered distribution of natural numbers on the square root spiral"**

   ( →  This study can be found under my author-name in the arXiv – data bank )